\documentclass[10pt]{amsart}
\usepackage{latexsym,amscd,amssymb,epsfig,psfig,epic,eepic} 
\usepackage[dvips]{color}

\theoremstyle{plain}
  \newtheorem{theorem}{Theorem}[section]
  \newtheorem{proposition}[theorem]{Proposition}
  \newtheorem{lemma}[theorem]{Lemma}
  \newtheorem{corollary}[theorem]{Corollary}
  \newtheorem{conjecture}[theorem]{Conjecture}
\theoremstyle{definition}
  \newtheorem{definition}[theorem]{Definition}
  \newtheorem{example}[theorem]{Example}
  
  \newtheorem{question}[theorem]{Question}
 \theoremstyle{remark}
  \newtheorem{remark}[theorem]{Remark}

\numberwithin{equation}{section}

\def\integers{{\mathbb Z}}

\begin{document}

\title[Connectivity of $h$-complexes]
{Connectivity of $h$-Complexes}

\author{Patricia Hersh}
\address{Department of Mathematics\\
     University of Michigan\\
     525 East University Ave.\\
     Ann Arbor, MI 48109-1109}
\email{plhersh@umich.edu}

\thanks{The author was supported by an NSF postdoctoral research
fellowship.}

\begin{abstract}
This paper verifies a conjecture of Edelman
and Reiner regarding the homology of the $h$-complex of a Boolean
algebra.  A discrete Morse function with no 
low-dimensional critical cells is constructed, implying a lower bound
on connectivity.  This together with an Alexander duality result of 
Edelman and Reiner implies homology-vanishing also in high dimensions.
Finally, possible generalizations to certain classes of supersolvable
lattices are suggested.
\end{abstract}

\maketitle

\section{Introduction.}

If a simplicial complex $\Delta $
has a shelling in which the unique minimal faces from the shelling steps
form a subcomplex of $\Delta $, then Edelman and Reiner call this
subcomplex the {\it $h$-complex} of $\Delta $ with respect to this
shelling.  They refer to such a shelling as an {\it $H$-shelling}.  Edelman
and Reiner introduced and studied $h$-complexes in [3].
One motivation
for $h$-complexes is that the $f$-vector of an $h$-complex is
the $h$-vector of the original complex and the Euler characteristic of
an $h$-complex is the Charney-Davis quantity of the original complex 
(cf. [13], [10], [11]). 

Following [4], let $\Delta_n$ denote the $h$-complex which results from
the standard shelling for the order complex of a truncated Boolean
algebra $B_n - \{ \hat{0},\hat{1} \}$.
Edelman and Reiner conjecture in [4] the following:

\begin{conjecture}[Edelman-Reiner]
$\tilde{H}_i (\Delta_n, \integers )$ is nonzero if and only if 
$ (3i + 5)/2 \le n \le 3i+4$.
\end{conjecture}
This is equivalent to saying that the reduced homology in
dimension $i$ is nonzero if and only if 
$$\frac{n-4}{3} \le i \le \frac{2n-5}{3}.$$
Our main result will be a proof of this conjecture.
Afterwards we suggest possible generalizations.  

Recall that a simplicial complex $\Delta $
is {\it pure} if all maximal faces are
equidimensional; these maximal faces are called its {\it facets}.  
A pure simplicial complex is {\it shellable} if there is a 
total order $F_1,\dots ,F_k$ on its facets with the following
property: for each $1\le j\le k$, there is a unique face $\sigma_j$
contained in $F_j$ which is minimal among all faces contained in $F_j$ but
not in any earlier facets.  
We refer to the faces $\sigma_1,\dots ,\sigma_k$ as the {\it minimal faces}
of the shelling.  For a shellable complex $\Delta $ of dimension $d$, 
the $h$-vector of $\Delta $ has
coordinates $(h_{-1},\dots ,h_d)$, with $h_i$ counting the number of 
facets $F_j$ for which the minimal face $\sigma_j$ is $i$-dimensional.

Our interest will be in the Boolean algebra $B_n$, namely the partial order
on subsets of $\{ 1,\dots ,n\} $ by inclusion.  Let $\hat{B_n}$ denote 
the truncated Boolean algebra $B_n - \{ \hat{0} ,\hat{1} \}$ consisting
of all subsets except the empty set and the full set.
Denote by $\Delta (P)$ the {\it order complex} of a poset $P$, i.e. the 
simplicial complex whose faces are the chains of comparable poset elements.
It is well-known that $\Delta (\hat{B_n})$ has a (lexicographic) shelling 
by labelling saturated chains with permutations in 
$S_n$ recording the order in which elements of $\{ 1,\dots ,n\} $ are
successively inserted, and then ordering facets in $\Delta (\hat{B_n})$
by the lexicographic order on the permutations labelling the saturated
chains.  The minimal faces in this 
shelling are comprised of the ranks at which the permutations have descents.
It is not hard to check that these minimal faces form a subcomplex,
denoted $\Delta_n$, of $\Delta (\hat{B_n})$.  

Reiner observed that the reduced homology $\tilde{H}_i(\Delta_n, \integers)$ 
is nonzero for $\frac{n-4}{3}\le i\le \frac{2n-5}{3}$ (personal
communication).  
A proof of his result is provided in Section ~\ref{non-vanish-section}.
In light of Reiner's observation, 
it will suffice to show that the homology vanishes in the remaining
dimensions.  We will use discrete Morse theory in Sections ~\ref{match-sec}
and ~\ref{acyclic-sec} to show this for 
dimensions below $\frac{n-4}{3}$.  
Then we use Theorem 4.14 (an Alexander duality result) from [4]
to deduce homology vanishing for top dimensions.  
Theorem 4.14 of [4] is as follows (see [4] for definitions):

\begin{theorem}[Edelman-Reiner]\label{Alex-dual}
Let $\omega $ be an $H$-shelling of a simplicial $d$-sphere $\Sigma $, and 
$\alpha $ a simplicial involution on $\Sigma $ which reverses the 
restriction map.  Denote by $\Delta^{(h)} (\omega )$ 
the $h$-complex given by 
$\omega $.  Then there is an isomorphism 
$$ \tilde{H}^i (\Delta^{(h)} (\omega ),\integers ) 
\rightarrow \tilde{H}_{d-1-i} (\Delta^{(h)} (\omega ), \integers).$$
\end{theorem}

Since the order complex of the truncated
Boolean algebra is the first barycentric subdivision of the boundary of
a simplex, it is a triangulation of a 
sphere.  Its standard shelling is an $H$-shelling, so
the above theorem applies to its $h$-complex.  

We will also give a combinatorial
proof that there is a dual discrete Morse function for 
the $h$-complex of a Boolean algebra in Section ~\ref{dual-sec}, 
yielding a second proof that 
its high-dimensional homology vanishes.  Our motivation for
this alternate proof is that it has the potential to generalize to 
situations where the Alexander duality result of [4] would not apply (e.g.
to the poset of subspaces of a finite vector space).  This dual Morse 
function might also be helpful for another question of [4], that of finding
a combinatorial 
explanation for the symmetry of the Betti numbers which results from 
Theorem ~\ref{Alex-dual} above.  Finally, Section ~\ref{open-sec} 
discusses possible generalizations from the Boolean algebra to other 
supersolvable lattices.

Before turning to the details, we quickly review the bare essentials 
from Forman's discrete Morse theory (see [5]) and Chari's combinatorial 
reformulation (see [3]).  See [2] for more background on topological
combinatorics and see [13] for background on $f$-vectors, $h$-vectors 
and the Charney-Davis conjecture.

\begin{definition}
A matching on the face poset $F(\Delta )$ of a simplicial complex 
$\Delta $ is {\it acyclic} if orienting 
matching edges upward and all other edges downward yields an acyclic
directed graph.  (Recall that $F(\Delta )$ is the partial order on 
faces by inclusion.)
\end{definition}
 
Any acyclic matching on $F(\Delta )$ gives rise to a 
discrete Morse function on $\Delta $ whose {\it critical cells} are
the faces left unmatched by the acyclic matching.  The number of 
critical cells of various dimensions
in a discrete Morse function give bounds on the Betti 
numbers as follows.  For each $i$, $\beta_i\le m_i$, 
where $m_i$ is the number of $i$-dimensional critical cells.

For simplicity, we will work exclusively with acyclic matchings rather than
the corresponding discrete Morse functions.  
Forman proved that $\Delta $ a discrete Morse function on a $d$-dimensional
CW complex $\Delta $ with Morse numbers $m_0,m_1,\dots ,m_d$ implies 
that $\Delta $ is homotopy equivalent
to a CW complex which has $m_i$ cells of dimension $i$ for each $i$.
We will specifically use the fact that a 
discrete Morse function on a complex
$\Delta $ with $m_i=0$ for $i$ less than a fixed
$j$ implies that the $\Delta $ is 
$(j-1)$-connected.  

\section{The $h$-complex of a truncated Boolean algebra}\label{setup-sec}

This section gives more detail about the standard shelling for the 
Boolean algebra $B_n$ in order to set up notation that we
will need for the acyclic matching in later sections.
The elements of $B_n$ are the 
subsets of $[n] := \{ 1,\dots ,n\} $.  $B_n$ has covering relations
$T\prec S$ for each $S= T\cup \{ i\} $ and $i\in [n]\setminus T$.
Label each covering relation $T\prec T\cup \{ i\} $
by the label $i$.  Each saturated 
chain is then labelled by 
the sequence of labels on its covering relations, i.e. by a permutation
in $S_n$ written in one-line notation.  Notice that each element of $S_n$
labels a single saturated chain, allowing us to refer to permutations
and saturated chains interchangeably.  Ordering these label
sequences lexicographically gives a shelling order on facets of
$\Delta (\hat{B_n})$.

Notice that the minimal face for a permutation $\pi $ 
consists of the chain 
supported at those ranks where $\pi $ has descents.  
For example, the minimal face for $\pi = 132654$ is the chain
$\{ 1,3\} < \{ 1,3,2,6\} < \{ 1,3,2,6,5\} $ which consists of the 
ranks of the descents $32, 65$ and $54$.  We can easily recover a saturated
chain from its minimal shelling face, so we also refer to minimal faces 
interchangeably with permutations and saturated chains.  It is
immediate from this description of minimal faces
that the Charney-Davis quantity is the alternating sum $A_n$ of the
Eulerian numbers $A_{n,k}$.  The exponential generating function 
$\sum_{n\ge 0} A_n \frac{x^n}{n!}$ is well-known to equal $-\tanh (x)$
(see [4, p. 52]).  Nonetheless, the
homology of $\Delta_n$ will turn out to live 
in many different dimensions.

To simplify notation later, we add an 
initial letter $a_0=0$ and a final letter $a_{n+1}=n+1$ to each permutation
$a_1\cdots a_n$.  
We refer to the permutation position between $a_i$ and $a_{i+1}$ as rank 
$i+1$, reflecting the fact that we have adjoined $a_0$ between ranks 0 and 1.
Depict the minimal shelling face 
$$\{ a_{1,1},\dots ,a_{1,i_1}\} < \{ a_{1,1},\dots ,a_{1,i_1},
a_{2,1},\dots ,a_{2,i_2} \} < \cdots < \{ a_{1,1} ,\dots ,a_{j,i_j}  \} $$ 
for a permutation
$\pi = a_0 a_1 \cdots a_{n+1}$ which has descents at ranks 
$i_1,i_1+i_2, \dots ,i_1+ \cdots + i_{j-1}$
by an ordered collection of blocks 
$$a_{1,1}\cdots a_{1,i_1} | a_{2,1}\cdots a_{2,i_2}| \cdots | a_{j,1}\cdots
a_{j,i_j} .$$ 
By convention, 
order elements within each block in increasing order.
Sometimes we refer to these blocks as {\it intervals}. 
Notice that this minimal face has dimension $j-2$.

We will call the separators between the blocks {\it bars}.
When we remove a bar and merge two consecutive blocks, we 
sort the two blocks so the new permutation is increasing on 
the merged block.  When we speak of inversions between two consecutive
blocks, we mean inversions in the permutation obtained by removing
the separating bar without sorting the blocks.

\section{Non-vanishing homology}\label{non-vanish-section}

This section shows for each integer $n$ with $\frac{3k+5}{2}
\le n\le 3k+4$ that the
homology group $\tilde{H}_k(\Delta_n,\integers )$ is nonzero.  
The approach for $2k+3\le n\le 3k+4$
is to exhibit a cycle which cannot be a boundary, by virtue of 
having a free face.  Theorem ~\ref{Alex-dual} 
then gives us that $\tilde{H}_k(\Delta_n , \integers )\ne 0$ 
for $\lceil \frac{3k+5}{2} \rceil \le n \le 2k+3$.  

\begin{definition}
A {\it free face} of dimension $k$ in $\Delta_n$ is a face which 
is not in the boundary of any $(k+1)$-dimensional face.  Thus, it
is a maximal face in $\Delta_n$.  
\end{definition}

First notice that $\Delta_n$ is not pure, so there will
be free faces in various dimensions.  Each cycle we construct in this 
section will contain at least one free face,
making it impossible for the cycle to be a boundary.

\begin{example}\label{cycle-extreme}
For $k=1,n=3k+4$, the minimal shelling face 
$13 | 246 | 57$ is a free face in $\Delta_7$.  
To see this, notice that any $2$-face 
containing this edge must permute elements within 
one or more of the blocks $13, 246, 57$ in a way 
that maintains the descents at ranks $2$ and $5$, and 
creates one new descent.  However, swapping $1$
and $3$ cannot avoid turning rank 2 into an ascent,
and likewise the descents at ranks 2 and 5 force
the labels $2,6$ and $5$ into their current positions,
making a $2$-face containing $13|246|57$ impossible.
Finally, notice that $13|246|57$ appears in the cycle
$$z = 13|246|57 - 13|26|457 + 3|126|457 - 3|1246|57.$$
\end{example}

The remainder of this section generalizes this to all $n,k$ satisfying
$2(k+1)+1 \le n\le 3(k+1) + 1$, by giving constructions in the two 
extreme cases, then showing how to combine them to yield the desired 
range.  Notice first that the cycle $z$ in 
Example ~\ref{cycle-extreme} may be viewed as a sum over
permutations of the form
$(12)^{e_1} (45)^{e_2}$ that act on positions, each applied
to the free face $C$.  Each permutation is multiplied by its
sign, to ensure that we get a cycle.  That is, 
$$ z = \sum_{\pi \in \langle (12), (45)\rangle } \mathrm{sgn} (\pi ) 
\pi (13|246|57) .$$  Note that the minimal shelling face for a 
permutation appearing in $z$ includes either rank 1 or 2, but not 
both, depending on whether $13$ appears in decreasing or increasing
order, and likewise includes either rank 4 or 5, depending on the 
order of 4 and 6.  Thus, the cycle is an alternating 
sum of $2^{k+1}$ faces, each of dimension $k$, chosen so 
that each $(k-1)$-face appearing in any of these $k$-faces
will occur in exactly two of them which have 
opposite signs.  This will ensure $\partial (z)=0$, as needed
for a cycle.

More generally, for $n=3k+4$ we use the free face
$$F = 13|246|\cdots | 3i+2,3i+4,3i+6  |\cdots |3k+2,3k+4 $$ in which 
the $i$-th block has elements $3(i-2)+2,3(i-2)+4,3(i-2)+6$ for each $i$  
strictly between $1$ and $k$.  A cycle $z$ is obtained by choosing
the order of the last two elements of each of the first $k+1$
blocks, i.e. for every block except the very last one.  Each of 
these pairs of block elements determines the location of one of 
the descents.  Thus, 
$$z = \sum_{\pi \in \langle (12),(45),(78),\dots ,(3k+3,3k+4)\rangle }
\mathrm{sgn} (\pi ) \pi (F)$$ with permutations $\pi $ acting on positions.

At the other extreme, for $n=2k+3$, one obtains a
free face $1,n|2,n-1|\dots |(n-1)/2,(n+3)/2|(n+1)/2$.  A
cycle again results from choosing the relative order of
$i,n-i+1$ for $1\le i \le k+1$, i.e. for pairs in each block except the
last one.  For $n$ satisfying $2(k+1)+1 \le n \le 3(k+1) + 1$,
we combine the two constructions above to obtain a free face
from the following permutation in $S_n$.  For
$n= 2j + 3(k+1-j) + 1$, begin the permutation with
$ 1,n|2,n-1|\cdots | j-1,n-(j-1)+1 $.  Appended to this is the following
permutation in $S_{[j,n-j+1]}$:
$$j,j+2| j+1,j+3,j+5 | \cdots | n-j-4,n-j-2,n-j| n-j-1,n-j+1 .$$
Again, we show this belongs to a cycle with $2^{k+1}$ faces by summing
over elements of a group of size $2^{k+1}$ each multiplied by its sign.
That is, for each of the first $k+1$ blocks, choose whether or not to
swap the order of the last two letters.
By similar reasoning to the above example, one obtains:

\begin{theorem}[Reiner]
For each $2(k+1) + 1 \le n \le 3(k+1) + 1$, $\tilde{H}_k(\Delta_n)\ne 0$.
Furthermore, by Theorem ~\ref{Alex-dual}, this implies
$\tilde{H}_k(\Delta_n)\ne 0$ for $\lceil \frac{3k+5}{2} 
\rceil \le n \le 2k+3$.
\end{theorem}

\section{A matching on $\Delta_n $}\label{match-sec}

This section provides a matching on faces in $\Delta_n$ which will be 
shown to be acyclic in the next section.  
In contrast to most acyclic matchings in the literature, our
matching is fairly easy to describe, but the proof of its acyclicity
is much more intricate than usual.  

Faces will be matched greedily
based on their lowest interval which takes a certain form, described 
below.  First we will need some notation.
Denote by $I^{above}$ the interval immediately above an interval 
$I$ when such an interval exists, and likewise denote by $I_{below}$ 
the interval immediately below $I$.  Let $S_C(I)$ be the size of 
the maximal set $S$ of consecutive
intervals $J_1,\dots ,J_s$ immediately above 
$I$ such that (1) $|J_1|=\cdots = |J_s| =2$, and (2)
the only inversions among the blocks
$I,J_1,\cdots ,J_s$ 
are the $s$ descents separating the blocks.  We call 
$J_1,\dots ,J_s$ the $J$-invervals or $J$-blocks of $I$.  
For example, in $0,1,2,3,6|5,8|7,9|4,10$ the block $I= 0,1,2,3,7$ has
$S_C(I)=2$ and has $J$-blocks $5,8$ and $7,9$ but not $4,10$.

\begin{definition}
An interval $I$ is {\it matchable} if it has any of the following forms:
\begin{enumerate}
\item
$|I| = 1$, $|I^{above}|$ is odd of size at least 3, and there is only
one inversion between $I$ and $I^{above}$.
\item
$|I|$ is even, $|I|\ge 4$, $I_{below}$ exists, and there are inversions 
between the largest element of $I_{below}$ and 
both of the two smallest elements
of $I$.
\item
$|I|\ge 4$, $S_C(I)$ is even, and $I$ is not also of type 2.
\item
$|I|\ge 2$, $S_C(I)$ is odd, there is only one inversion between $I$
and $I^{above}$, and the block obtained by merging $I$ with
$I^{above}$ is not matchable of type 2.
\end{enumerate}
\end{definition}

In an effort to make our proofs more readable, let us call the four
types of matchable intervals above (1) 1-split, (2) 1-merged, 
(3) 2-merged, and (4) 2-split, respectively, reflecting the fact that
a block of size 1 or 2 is split off from another block or merged with
it.  Notice that $0$ and $n+1$ are permanently fixed in the first and
last positions, so the matching may not
insert bars at ranks $1$ and $n+1$; the requirement for 1-merged blocks
$I$ that $I_{below}$ exists will take care of this. 

When we need
to keep track of the fact that we are viewing $I$ as an interval in a 
chain $C$, then we will sometimes denote $I$ as $I(C)$.
If the first matchable interval in a chain is at rank $r$, then we
match it with another chain whose first matchable interval is also 
at rank $r$, as follows.  

\begin{definition} A chain $C$ with lowest rank matchable interval $I(C)$
at rank $r$ is {\it matched} with a chain $D$ if $D$ differs from $C$ by a 
single inversion and:
\begin{itemize}
\item
$I(C)$ has type 1; $D$ is obtained from $C$ by merging $I$ with 
$I^{above}$.
\item
$I(C)$ has type 2; $D$ is obtained from $C$ by splitting $I$ of size 
$m$ into blocks of size $1,m-1$ (where we list the block at higher
ranks second).
\item
$I(C)$ has type 3; $D$ is obtained from $C$ by splitting $I$ of size 
$m$ into blocks of size $m-2,2$.
\item
$I(C)$ has type 4; $D$ is obtained from $C$ by merging $I$ with 
$I^{above}$.
\end{itemize}
\end{definition}

To prove the above matching is well-defined, we first show that if $C$ has
lowest matchable interval at rank $r$, then its partner $D$ also has
lowest matchable interval at rank $r$.  

\begin{theorem}
If the lowest matchable 
interval in a chain $C$ is at rank $r$, then the chain $D$ with which $C$
is matched also has no matchable intervals below rank $r$.  
\end{theorem}

\begin{proof}
Suppose the lowest matchable
interval $I$ in $C$ is 1-split.  Then
$|I(C)|=1$, and $|I(D)|$ is even with size at least 4.
Since $C$ and $D$ agree below rank 
$r$ and all intervals of size at least 4 are matchable, there 
cannot be any intervals of size 4 or larger in $D$ below 
rank $r$.  Hence, $D$ has no 1-merged or 2-merged matchable
intervals below 
rank $r$.  Neither $I(C)$ nor $I(D)$ has 
size 2, so neither can be a $J$-interval for any lower intervals.
$D$ cannot have a 2-split matchable interval $I'$ at rank $r'<r$ without
$C$ also having such an interval: $C$ and $D$ agree below rank $r$, and 
$I'(D)$ cannot
have any $J$-intervals at or above rank $r$, so $S_D(I') = S_C(I')$ .
Finally, suppose
$D$ had a 1-split 
matchable interval $I'$ at rank $r'<r$.  Then $I'(C)$
would also be matchable, except perhaps for $r'=r-1$.
But then $D$ would need an interval 
of odd size at rank $r$, but $|I(D)|$ is even.  Hence, 
$D$ has no matchable intervals below rank $r$.  The case where $I(C)$ is
1-merged is similar with the roles of $C$ and $D$ reversed, 
so we omit the argument.

Now suppose $I(C)$ is 2-merged or 2-split.
Once again $C$ and $D$ agree below 
rank $r$, and all intervals of size at least four are matchable; thus,
we only need to consider the possibility that $D$ has a matchable 
interval $I'$ at rank $r' < r$ with $I'$ that is 1-split or 2-split.
If $I'(D)$ is 1-split matchable,
then as before $I'$ must occur at rank $r-1$.  Then $|I(D)|=m$ is 
odd with $m\ge 3$, which means $|I(C)| = m\pm 2$ is also odd.  Furthermore,
$|I(C)|\ge 2$ which means it will also have size at least 3, since 
$|I(C)|$ is odd.  Furthermore, $I'(C)=I'(D)$ and $I'_{below}(C) = 
I'_{below}(D)$, so $I'(C)$ would 
also be 1-split matchable, a contradiction.

Now suppose $I'(D)$ is 2-split matchable.
Then we would need $S_D(I')$ odd and $S_C(I')$ even, so in particular they
are not equal.  Since
$C$ and $D$ agree below rank $r$, this means that either $I'(C)$ or
$I'(D)$ must have one or more $J$-intervals at or above rank $r$.
Hence, $C$ or $D$ must have a block of size 2 at rank $r$, while the
other must then have a block of size 4 at rank $r$.
Let us assume $|I(C)|=4$, which means $S_C(I)$ is even.  The other case
is similar.  

Since 
$S_C(I')$ is even and $C$ does not have a block of size 2 at rank $r$, both
$I'(C)$ and $I'(D)$ must have an even number of $J$-blocks below rank $r$.
Thus, $I'(D)$ needs an odd number of $J$-blocks above rank $r$.  However,
$S_D(I)$ is odd, implying $I(D)$ together with its $J$-blocks comprise an 
even number of  
prospective $J$-blocks for $I'(D)$ above rank $r$.  This means that
not all of the $J$-blocks for $I$ are also $J$-blocks for $I'$, so there
must be at least one extra inversion among these potential $J$-blocks.
In particular, either the second smallest label
above rank $r$ must be smaller than the label just below rank $r$, or else
the smallest label above rank $r$ must be smaller than the second smallest
label below rank $r$.  We can eliminate the latter possibility, since 
$I'(D)$ has at least one $J$-block above rank $r$.
Hence, $I(C)$ has two labels that are smaller than the 
largest element of $I_{below}(C)$, and $|I(C)|$ is even of size at least 4.
This means that $I(C)$ is 1-merged matchable instead of 2-merged matchable, 
a contradiction.
\end{proof}

\begin{corollary}
The matching is well-defined.
\end{corollary}

\begin{proof}
It suffices now to check that the matching rules for 1-split and 1-merged
matchable intervals are
inverses to each other, and likewise for 2-merged and 2-split intervals.  
This is easy, and is left to the reader.  
\end{proof}

\section{Acyclicity of $\Delta_n $ matching}\label{acyclic-sec}

Now we turn to the task of proving the matching is 
acyclic, and hence comes from a discrete Morse function.  Unlike many
acyclicity proofs in the literature, we are not aware of any function which
is decreasing along directed paths, so our acyclicity proof will take
another approach.  

\begin{lemma}\label{length-preserve}
If the matching had a directed cycle $C$, then each downward 
step in $C$ would eliminate a single inversion, i.e. would 
merge two blocks with only one inversion between them.
\end{lemma}

\begin{proof}
Each upward step increases permutation
length (i.e. number of inversions)
by exactly one, and each downward step decreases permutation length
by at least one.  Any cycle would have an equal number of upward and 
downward steps before revisiting its initial permutation, so down steps
must decrease length by exactly one, in order to restore the length of 
the original permutation.
\end{proof}

In light of Lemma ~\ref{length-preserve}, each edge traversed in a 
directed cycle may be viewed as an adjacent
transposition; an entire cycle would comprise a non-reduced
expression for the identity permutation.  It would be 
desirable to have a shorter, more elegant
proof of acyclicity than the one below, perhaps using properties
of non-reduced expressions for the identity permutation.  

Before proceeding with the proof, 
we list a few facts it will use repeatedly:
\begin{enumerate}
\item
By Lemma ~\ref{length-preserve}, 
downward 
steps merging two blocks are only permitted when the only inversion
between the blocks is the descent separating them.
\item
Since each upward step changes the lowest matchable interval from 2-merged 
or 1-merged to 2-split or 1-split, it must be immediately followed by a 
downward step which causes
the lowest matchable interval to again be 2-merged or 1-merged.  Otherwise
the downward step could not be followed by another upward step, as would be 
required in a cycle.  
\item
There are no matching
steps splitting a block of size $m$ into smaller blocks of size $m-1,1$.
\end{enumerate}
Within the proof, we refer to these facts as Observations 1, 2 
and 3.

One other key ingredient will 
be the idea behind the 0-1 Sorting Lemma from 
theoretical computer science (cf. [8]), that deals with the following
type of sorting procedure: an {\it oblivious comparison-exchange
sorting procedure} is an ordered 
list of comparisons to be performed, where two
elements are exchanged whenever they are compared and found to be out
of order; this is ``oblivious'' in that the choice of comparisons cannot
depend on the outcome of earlier comparisons.

\begin{lemma}[0-1 Sorting Lemma]
Any oblivious comparison-exchange sorting
algorithm which correctly sorts lists consisting exclusively
of 0's and 1's will correctly sort 
lists with arbitrary values.  
\end{lemma}

The idea is that to
sort numbers correctly, one must be sure for any particular value $a$ that
all numbers larger than $a$ are sorted to above all numbers smaller than
$a$, and so for any fixed $a$ one may treat the numbers larger than $a$ as
1's and those smaller than $a$ as 0's.  In our context,
we will have a particular label $a$ and it will be quite useful to keep
track of exactly which labels below it form inversions with it, and 
to disregard all other information about the relative
order of the values below $a$.

\begin{remark}
The proof below often speaks of rank, by which we mean rank in the original
poset $\hat{B_n}$, not in the face poset $F(\Delta_n)$ upon which we
construct a matching.
\end{remark}

Denote by $u_r$ the matching step which inserts a bar at rank $r$.  Denote
by $d_r$ the downward step deleting a bar from rank $r$ by applying an
adjacent transposition to replace a descent by an ascent.  

\begin{theorem}\label{acyclic}
The matching on $\Delta_n$ is acyclic.
\end{theorem}

\begin{proof}
Suppose there were a directed cycle $C$ in the directed graph obtained 
from the matching on $F(\Delta_n)$.  Consider the highest rank $t$
at which a bar is ever inserted, and let $u_t$ be a matching 
step inserting such a bar $B_t$ into a chain $C_0$ to obtain a partner 
chain $D_1$.  
\begin{figure}[h]
\begin{picture}(150,70)(105,0)
\put(10,10){\line(1,1){30}}
\put(60,40){\line(1,-1){30}}
\put(110,10){\line(1,1){30}}

\put(210,40){\line(1,-1){30}}
\put(260,10){\line(1,1){30}}
\put(310,40){\line(1,-1){30}}
\dashline{2}(160,25)(190,25)

\put(0,-5){$C_0$}
\put(45,50){$D_1$}
\put(95,-5){$C_1$}
\put(145,50){$D_2$}

\put(195,50){$D_{k-1}$}
\put(245,-5){$C_{k-1}$}
\put(295,50){$D_k$}
\put(345,-5){$C_k$}

\put(10,30){$u_t$}
\put(110,30){$u_{i_2}$}
\put(260,30){$u_{i_k}$}
\put(330,30){$d_t$}

\end{picture}
\caption{The cycle segment from $u_t$ to $d_t$}
\label{cycle-segment}
\end{figure}
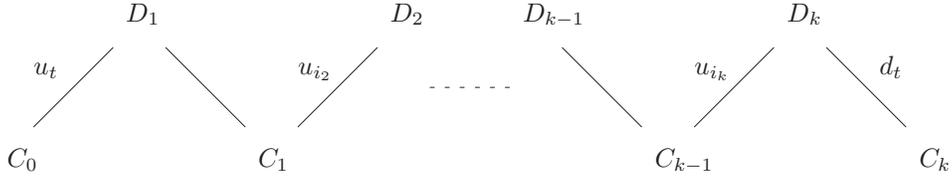
Let $u_{i_k}$ be the upward step immediately preceding the 
first occurence of $d_t$ after $u_t$.  Then $d_t$ deletes a bar at a 
strictly higher rank than $i_k$ (since $i_k \le t$, but we are assured
that $i_k\ne t$ since there is already a bar at rank $t$ just prior to 
$u_{i_k}$).  

Our proof
will focus on the segment of $C$ from just before $u_t$ until just after
$d_t$.  Let us establish some notation for the faces appearing in 
this segment.  $C$ must alternate between two consecutive face poset ranks 
$r$ and $r+1$, so denote this segment of $C$ by
$C_0 \rightarrow D_1 \rightarrow C_1 \rightarrow \cdots \rightarrow
C_{k-1} \rightarrow D_k \rightarrow C_k$.  That is, denote chains
at rank $r$ by $C_0,C_1,\dots , C_k$ and chains at rank $r+1$ by 
$D_1, \dots ,D_k$;
the $j$-th matching step in this segment, denoted $u_{i_j}$, 
takes $C_{j-1}$ to $D_j$.  We chose $k$ so that the first $d_t$ after
$u_t$ immediately follows $u_{i_k}$, so there are $k$
matching steps within the segment.  See Figure
~\ref{cycle-segment}.

Now we already observed that $d_t$ deletes a bar at a strictly higher
rank than where a bar was inserted by $u_{i_k}$.  However, since 
$u_{i_k}$ is a matching step, it must have changed the 
lowest matchable interval $I$ 
in $C_{k-1}$ from 1-merged to 1-split or from 2-merged
to 2-split.  Since $t>i_k$, $d_t$ cannot create a lower matchable interval
than $I(C_{k-1})$.  Thus, by Observation 2, 
$d_t$ must destroy the structure
which made $I(C_{k-1})$ 1-split or 2-split matchable. 
We will consider two cases, depending
on whether $u_{i_k}$ changes an interval $I$ either (a) from 1-merged
to 1-split
or (b) from 2-merged 
to 2-split.  Each case will lead to a contradiction, making
cycles impossible.  

Proof in case (a):  
Suppose $u_{i_k}$ changed an interval $I$ from 1-merged to 1-split.  Lemma
~\ref{lemma-A} 
will show $u_{i_k}$ splits the block $I_{t^-}$ 
immediately below $B_t$.  Thus,
the block created immediately above $B_t$ 
has odd cardinality.  
Thus, the block created by $u_t$ just above $B_t$ 
cannot have size 2, so $u_t$ must 
have changed an interval from 1-merged to 1-split.
Therefore, $C_0$
must have a bar at rank $t-1$.  However, $C_k$ cannot have a bar at rank
$t-1$, since $u_{i_k}$ split $I_{t^-}$ into 
two blocks, with the one at higher ranks having size at least 3.  Thus, 
there must have been an intermediate step $d_{t-1}$ at some point.  However,
Lemma ~\ref{lemma-B} show 
that steps $u_{t-1}$ are impossible throughout the cycle, a 
contradiction to our ever returning to $C_0$.  

Proof in case (b): 
Suppose alternatively that $u_{i_k}$ changes a
matchable interval $I(C_{k-1})$ 
from 2-merged
to 2-split.  This is immediately followed by $d_t$ which
deletes a bar strictly above rank $i_k$, but
by Observation 2, $d_t$ 
causes $I(D_k)$ no longer to be 2-split matchable.
To do this, $d_t$ must change the parity of 
$S_{D_k}(I)$.  Hence, $d_t$ must delete a bar separating a $J$-block
just below $B_t$ from a non-$J$-block $I_{t^+}$
immediately above $B_t$.  To avoid
being a $J$-block, $I_{t^+}$ must either (I) have size 
$m > 2$ or (II) have size $m=2$ and have an inversion 
with the $J$-blocks below it
other than the descent separating $I_{t^+}$ from $I_{t^-}$.

Case b(I): 
If $m >2$, then $u_t$ must have changed a matchable interval 
from 1-merged to 1-split.
Notice that there is a bar at 
rank $t-1$ just prior to $u_t$ but no bar at rank $t-1$ just after
$u_{i_k}$ (since there is a $J$-block immediately below 
rank $t$ just before $d_t$).  
This means we need a step $d_{t-1}$ in our cycle, but Lemma ~\ref{lemma-B}
again such a step.

Case b(II):
If $m=2$, there must be extra inversions preventing the block 
$I_{t^+}$ just above rank
$t$ from being a $J$-block for $I(D_k)$.  Lemma ~\ref{lemma-C} will show 
that the larger element $d$ in  $I_{t^+}$
cannot be inverted with any 
elements of the $J$-blocks of $I(D_k)$, 
so that the extra inversion must instead involve the smaller element
$a$ in the block above $B_t$.  Thus, 
$a$ must be smaller than the two largest 
labels within 
$J$-blocks of $I(D_k)$, i.e. the labels just 
below ranks $t$ and $t-2$ in $D_k$.  
Hence, $D_k$ has an inversion between the label $a$ just above rank $t$ 
and the label just below rank $t-2$.  However, Proposition ~\ref{prop-G}
will show we cannot get from
the chain $C_0$ whose lowest matchable interval is 2-merged
to the situation at $D_k$ where the letter $a$ just above $B_t$ is
inverted with the letter just below rank $t-2$.  This will
complete our proof.
\end{proof}

\begin{lemma}\label{lemma-A}
Let $u_t$ insert bar $B_t$ in the highest position a bar is ever inserted
within a directed cycle.  If the step $u_{i_k}$ immediately preceding the
next $d_t$ after $u_t$ changes an interval from 1-merged to 1-split, then
the block above $B_t$ has odd size.
\end{lemma}

\begin{proof}
This is because $d_t$ must cause $I$ no longer
to be 1-split, by Observation 2, and the
only way to do this is for $d_t$ to delete a bar immediately above 
$I^{above}(D_k)$
so as to change the parity of $I^{above}(D_k)$ from odd
to even; this can only be done by merging $I^{above}(D_k)$ 
with an odd block immediately above it.  
This odd block immediately above $B_t$ is left unchanged by all steps 
between $D_1$ and $D_k$,
implying that $u_t$ inserted a bar $B_t$ making
the block immediately above $B_t$ odd.  
\end{proof}

\begin{lemma}\label{lemma-B}
If the highest insertion $u_t$ in a directed cycle changes an interval
from 1-merged to 1-split, then steps $u_{t-1}$ are impossible in the 
cycle.
\end{lemma}

\begin{proof}
The step $u_t$ split a block of even size $n$ into blocks 
of size $1,n-1$.  By Observation 3, we cannot have
an upward step $u_{t-1}$ while a bar is present at rank $t$, since no
matching steps split an interval of size $m$ into smaller intervals of
size $m-1,1$, where we list the higher interval second.  On the other
hand, when $B_t$ is not present, then inserting $u_{t-1}$ would create a 
block of size $n-2$ above $B_{t-1}$.  However, $n-2$ must be even of size
at least 4, since $n-1$ was odd of size at least 3.  There are no such 
matching steps, so $u_{t-1}$ is impossible.  
\end{proof}

\begin{lemma}\label{lemma-C}
In case b(II), the larger element $d$ in the block above $B_t$ cannot
be inverted with any elements of the $J$-blocks of $I(D_k)$.
\end{lemma}

\begin{proof}
The point will be to show
that the label $c$ just below rank $t$ in $D_1$ must still be in this
position in $D_k$.  But then we know that $u_t$ only increased the 
permutation length by exactly one, so that $c<d$ since $c$ was in the 
same block with $d$ just prior to $u_t$.  Since $c$ must be larger than all
other elements of the $J$-blocks of $I$, $d$ must also be larger than all
of them.

To show that $c$ is still at the position just below $B_t$ 
in $D_k$, we will show that there could not have been a step $u_{t-1}$
or $d_{t-1}$ between $D_1$ and $D_k$.
By Observation 3 we could not have inserted a bar
at rank $t-1$ in this interval, because a bar was present at rank $t$ the
entire time.  On the other hand, $u_t$ must have changed a matchable interval
from 2-merged to 2-split, since $m=2$ and in particular is even; this implies
that $D_1$ does not have a bar at rank $t-1$ available to 
be deleted. 
\end{proof}

\begin{proposition}\label{prop-G}
It is impossible in case b(II) above to have a directed path from
the face $C_0$ in which the lowest matchable interval must have been  
2-merged to the face $D_k$ where the letter $a$ just above $B_t$ is
inverted with the letter just below rank $t-2$.
\end{proposition}

\begin{proof}
Lemma ~\ref{size-four} will show that 
$u_t$ must have split a block of size 4 into
blocks of size $2,2$.
Consider the cycle element $C_0$ just prior to $u_t$.  Denote by $K$ the
block just below the bar $B_{t-2}$ in $D_1$.  (Note that $D_k$ has a bar
at rank $t-2$, because 
$u_t$ split a block of size 4 into blocks of
size $2,2$.)  For $u_t$ to change a matchable interval from 2-merged to
2-split instead of from 1-merged to 1-split, we need the largest element
of $K$ to be smaller than $a$.  We also know that the element just above
rank $t-2$ in $D_1$ is smaller than $a$, since $u_t$ only increased the 
permutation length by one.  We consider two cases, depending
on whether (i) $|K|=1$ or (ii) $|K|\ge 2$.

Case (i): $|K|=1$, so $C_0$ has a bar at rank $t-3$ as well as 
rank $t-2$.
$D_k$ does not have a bar at rank $t-3$, since $u_{i_k}$ matches a 
2-merged matchable block, and then $d_t$ causes this block to no longer
be 2-split matchable, which means all blocks between the bar inserted 
by $u_{i_k}$ and $B_t$ have size 2.  Thus, 
there must be an intermediate step $d_{t-3}$ prior to 
$u_{i_k}$.  Next we 
show that the cycle can never restore the situation
of having bars at both ranks $t-3$ and $t-2$, which will give us 
a contradiction. First note that a bar cannot be inserted at rank 
$t-3$ while one is present
at rank $t-2$.  On the other hand, we cannot have a step $u_{t-2}$ while a 
bar is present at rank $t-3$, by the following reasoning: such a 
step would change a matchable interval from 1-merged to 1-split, so it would
product an odd block of size at least 3 immediately above rank $t-2$; this
is impossible both when a bar is present at rank $t$ and when there is no
bar at rank $t$, since then the next lowest bar is at rank $t+2$.  Thus,
(i) is impossible since we cannot return to having bars at both
ranks $t-3$ and $t-2$.

Case (ii): $|K|\ge 2$ and
all elements of $K$ must be smaller than the label $a$ appearing
just above $B_t$.  
In the spirit of the 0-1 Sorting Lemma, we now
denote numbers below rank $t$
as 1's and 0's depending on whether they are 
larger or smaller than this fixed value $a$.  Regardless of the 
actual values, any 1 with a 0 above it must
be a descent, while any 0 with a 1 above it is an ascent.  
Immediately after
$u_{i_k}$ we need there to be $1$'s just below ranks $t$ and $t-2$.
However, in $D_1$ we know that the block between ranks $t-2$ and rank $t$
consists of one 1 and one 0, while the block $K$ just below this contains
only 0's.
Hence, the 1 just below rank $t-2$ in $D_k$ must have moved 
upward from below the $K$ block.  Finally, Lemma ~\ref{lemma-E} will 
use the idea of the 0-1 Sorting Lemma to show that this is impossible,
again precluding a cycle.
\end{proof}

\begin{lemma}\label{lemma-E}
It is impossible in Case (ii) of Proposition ~\ref{prop-G}
for a directed path to proceed from
the face $D_1$ to the face $D_k$.  That is, we cannot shift a label 
which is larger than $a$ upward from below the $K$ block to just below
rank $t-2$.
\end{lemma}

\begin{proof}
In $D_1$, the highest 1 below rank $t-2$ must be below rank $t-4$, because it
must be strictly below the block $K$ in order for $I(C_0)$ to avoid being
1-merged, 
\begin{figure}[h]
\begin{picture}(60,220)(0,0)

\put(0,40){\line(1,0){20}}
\put(0,80){\line(1,0){20}}
\put(0,150){\line(1,0){20}}
\put(0,190){\line(1,0){20}}

\dashline{2}(10,-5)(10,20)
\dashline{2}(10,45)(10,75)
\dashline{2}(10,85)(10,110)

\put(8,25){$1$}
\put(8,117){$0$}
\put(8,135){$0$}
\put(8,157){$0$}
\put(8,175){$1$}
\put(8,197){$a$}
\put(8,215){$d$}
\put(25,185){$B_t$}
\put(25,145){$B_{t-2}$}
\put(40,108){$K$ of size at least 2}
\put(40,25){Highest $1$ below $K$}

\end{picture}
\caption{1's and 0's below $a$}
\label{high-star}
\end{figure}
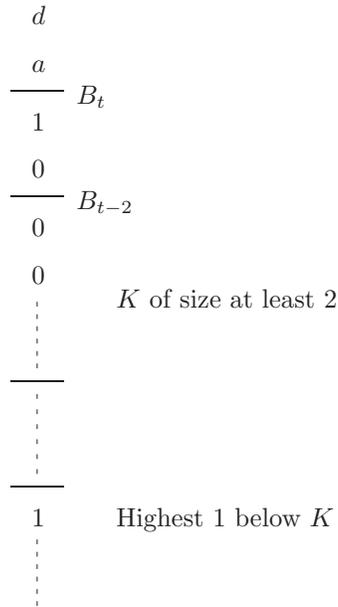
and we already showed $K$ has size
at least 2.  See Figure ~\ref{high-star}.

Furthermore, the only 1 in 
the interval between ranks $t-2$ and $t$ in $D_1$ is the one just below
rank $t$ that never moves.  Thus, we must eventually move
a 1 upward from below rank $t-4$ to just below rank $t-2$.  Just before 
moving a 
1 upward to just below rank $t-3$, we must have a bar at rank $t-3$, since
otherwise the step $d_{t-4}$ would eliminate more than one inversion, since
the 1 will be larger than all the 0's in the block above it.
However, there is no bar $B_{t-3}$ in $D_1$,
since $|K|\ge 2$.  Thus, we need a 
step $u_{t-3}$, and this can only happen when there is no bar at rank $t-2$,
by Observation 3.  Once we have a 
bar at rank $t-3$ with a 1 immediately below it, there will henceforth
be a 1 at this position until there is a step $d_{t-3}$, since bars
cannot be inserted at rank $t-4$ while a bar is present at rank $t-3$.
However, we cannot have $d_{t-3}$ until after $u_{t-2}$, since again 
$d_{t-3}$ would otherwise eliminate more than one inversion. 
Finally, it is not possible to have an upward step $u_{t-2}$ with a bar
present at rank $t-3$, again by a parity argument: $u_{t-2}$ would
need to change a matchable interval from 1-merged to 1-split, 
meaning we would
need a block of odd size at least 3 immediately above $B_{t-2}$, which is 
not possible since there is a bar at rank $t$.  
Thus, a directed cycle cannot get from the situation just
after $u_t$ to the situation needed just prior to $d_t$, contradicting 
there being a cycle.
\end{proof}

\begin{lemma}\label{size-four}
If $u_{i_k}$ and $u_t$ both change matchable intervals from 2-merged to 
2-split, then $u_t$ specifically must split a block of size 4 into blocks of 
size $2,2$.
\end{lemma}

\begin{proof}
Let $m$ be the size of the lowest matchable interval $I_0$ in $C_0$.
The idea of this lemma is that if $m-2 > 2$, then the second largest
element $b$
in the block below the bar inserted by $u_t$ is smaller than both labels
$c, e$ in the block above this same bar.  We cannot insert a 
bar at rank $t-1$ while a bar is present at rank $t$, and the label $b$ will
not move until we insert a bar at rank $t-2$.  Until such a bar is inserted,
all labels below $b$ in the block containing $b$ must be 
smaller than both $c$ and $e$, because they are smaller than $b$.

We must eventually insert a bar at rank $t-2$, since such a bar is 
present just after $u_{i_k}$.  However, we cannot have a matching step
inserting such a bar, under our $m-2>2$ assumption, since such a step
inserting a bar into a block $I$ would require $S_C(I)$ to be even, 
and we can show that $S_C(I)$ must be odd, as follows.  We have that
$S_{C_0}(I)$ was even, and we check next that 
$S_C(I) = S_{C_0}(I_0) + 1$.  

An interval above rank $t+2$ is a 
$J$-interval for $I$ if and only if it is a $J$-interval for $I_0$, 
since in either case the only allowable inversion between such 
intervals and elements in $I$ or $I_0$ is a single inversion with $e$.
Thus, $S_C(I) = S_{C_0}(I_0)+1$, so 
$S_C(I)$ cannot be even, a contradiction to $m-2$ being larger than 2.
\end{proof}

\section{Vanishing homology and a dual Morse function}\label{dual-sec}

\begin{theorem}
The $h$-complex $\Delta_n $ has a discrete Morse function with 
$m_i=0 $ for $3i + 4 < n $, so $\Delta_n $ is  $\lfloor \frac{n-5}{3} 
\rfloor $-connected.
\end{theorem}

\begin{proof}
Theorem ~\ref{acyclic} proves that our matching is acyclic, and hence
gives rise to a discrete Morse function whose critical cells are the 
unmatched face poset elements.  Since any interval of size at least four
is matchable, critical cells must have block sizes
$i_1,i_2,\dots ,i_j\le 3$ for
$i_1 + \cdots + i_j = n+2$ (recalling that we adjoined 
$a_0$ and $a_{n+1}$, increasing permutation lengths to $n+2$ letters).  
Hence, $3j \ge n+2$ for any unmatched face of dimension $j-2$, so
$m_j=0$ for $3j+4<n$, as desired.  
\end{proof}

Now we apply the Alexander duality of [4] result to deduce that there is also
no reduced homology in the necessary top dimensions.  Alternatively, this
may be verified by dualizing our matching construction, as follows.

\begin{theorem}\label{dual-morse}
$\Delta_n$ has a discrete Morse function with no critical cells of dimension
$i$ for $i>(2n-5)/3$, so $\tilde{H}_i (\Delta_n) = 0$ for $i>(2n-5)/3$.
\end{theorem}
\begin{proof}
We reverse the roles of ascents and descents in the original 
matching.  That is,
break any permutation into maximal blocks of decreasing labels, and put 
bars at the locations of all the ascents in the permutation.  Thus, bars are
at the ranks which are absent in the associated minimal shelling face.  We
may use the same matching construction as before, but with respect to 
this new choice of bars and blocks for each permutation.  Since matching 
steps inserting a bar will now eliminate exactly one inversion, more
specifically a descent, there is
no problem with having the bars at the missing ranks rather than
at the ranks present in a face.  Now all the arguments of the previous 
sections go through unchanged.  In conclusion, there are no critical
cells with four or more consecutive decreasing labels, 
implying there are no critical cells above dimension $(2n-5)/3$.
\end{proof}


\begin{question}
Is there a nice description of the permutations giving rise to 
critical cells?  Do some nice 
subset of these index a homology basis?  
Can we further collapse to this basis by gradient path 
reversal?  
\end{question}

\section{Possible generalizations}\label{open-sec}

Peter McNamara recently showed in [9] that supersolvability 
for a lattice of rank $n$ is equivalent to it
having an EL-labelling in which each edge is labelled by an integer
in $\{ 1,\dots,n\}$ in such a way that each saturated
chain is labelled with a permutation in $S_n$.  He calls such an 
EL-labelling an $S_n$ EL-labelling.  Richard Stanley previously provided an 
$S_n$ EL-labelling for every supersolvable lattice in [12].  It is 
shown in [4] that labellings known as 
SL-labellings (originally introduced in [1])
give $h$-shellings, and that supersolvable lattices 
have SL-labellings, namely their $S_n$ EL-labellings.  

\begin{question}  If $\Delta $ is the $h$-complex of a supersolvable
lattice of rank $n$ whose M\"obius function is nonzero on every interval,
then is $\Delta $ at least $\lfloor \frac{n-5}{3}\rfloor $-connected?
\end{question}

It seems plausible that $S_n$ EL-labellings might enable
one to generalize the discrete Morse function of previous sections to 
other supersolvable lattices.  The above 
M\"obius function requirement ensures that every interval has at least
one decreasing chain.  This seems essential to a matching in which 
all chains which include blocks of size 4 or larger are indeed matched.

\begin{remark}
The lattice of subspaces of a finite-dimensional vector space over a 
finite field is probably easier than the general question of any
supersolvable lattice with nowhere-zero M\"obius function.  Another
specific candidate would be the intersection lattice of any supersolvable 
arrangement.
\end{remark}

The following lemma from [7] seems likely to be helpful, in conjunction 
with a filtration by partially ordering Boolean algebras (e.g. apartments
in the poset of subspaces of a finite vector space).

\begin{lemma}[Cluster Lemma]\label{filter}
Let $\Delta $ be a regular $CW$ complex which decomposes into collections 
$\Delta_{\sigma }$ of cells indexed by the elements $\sigma $ in a partial 
order $P$ with unique minimal element $\hat{0}=\Delta_0$.  
Furthermore, assume
that this decomposition is as follows:
\begin{enumerate}
\item
$\Delta $ decomposes into the disjoint union $\cup_{\sigma\in P} 
\Delta_{\sigma }$, that is, 
each cell belongs to exactly one $\Delta_{\sigma }$ 
\item
For each $\sigma\in P$, $\cup_{\tau \le \sigma} 
\Delta_{\tau}$ is a subcomplex of $\Delta $
\end{enumerate}
For each $\sigma\in P$, let 
$M_{\sigma }$ be an acyclic matching on the subposet $F(\Delta |_{
\Delta_{\sigma }})$ of $F(\Delta )$ consisting of the cells in 
$\Delta_{\sigma }$.
Then $\cup_{\sigma\in P} M_{\sigma }$ is an acyclic 
matching on $F (\Delta )$.
\end{lemma}

Topologically, the order complex of
a supersolvable lattice with nowhere-zero M\"obius function
will consist of overlapping spheres, 
specifically overlapping type A Coxeter complexes.  We refer readers to
[6] for a potentially useful way of viewing those chains in
a Boolean algebra that do not belong to any earlier Boolean algebra as 
an intersection of half-spaces restricted to a sphere. 

\begin{remark}
The Alexander duality result of [4] will not apply to most supersolvable
lattices with nowhere-zero M\"obius function, since these will not in
general be spheres.  However, there could still be 
a dual discrete Morse function, similar to Theorem
~\ref{dual-morse}.
\end{remark}

\begin{question}
Is there a more general lower bound on connectivity for
$h$-complexes of SL-shellable posets 
whose M\"obius function is nonzero on every interval?
\end{question}

\section*{Acknowledgments}

The author thanks Vic Reiner and John Shareshian for 
helpful discussions.

\end{document}